\documentclass[11pt]{article}
\usepackage{amsmath,amssymb,amsthm,geometry,booktabs,longtable}
\newtheorem{theorem}{Theorem}[section]
\newtheorem{corollary}[theorem]{Corollary}
\newtheorem{proposition}[theorem]{Proposition}
\newtheorem{lemma}[theorem]{Lemma}
\newtheorem{conjecture}[theorem]{Conjecture}
\newtheorem{remark}[theorem]{Remark}
\title{Universal Implications Among $k$-Prime Graph Labelings}
\author{Jing Zhang\\
Department of Mathematics, Governors State University\\
University Park, Illinois 60484, USA\\
\texttt{jzhang@govst.edu}
\and
Yuan Li\\
Department of Mathematics, Winston--Salem State University\\
Winston--Salem, North Carolina 27110, USA\\
\texttt{liyu@wssu.edu}}
\date{}
\begin{document}
\maketitle
\begin{center}\small\textit{2020 Mathematics Subject Classification:} 05C78, 05C60, 11N05.\quad\textit{Keywords:} $k$-prime labeling, graph labeling, coprime graph, universal implication, prime switching, natural density.\end{center}
\begin{abstract}
For $k\ge1$, let $\mathcal P_k$ be the class of graphs admitting a bijective labeling by consecutive integers beginning at $k$, with coprime labels at adjacent vertices. We determine all universal inclusions between distinct such classes: $\mathcal P_i\subseteq\mathcal P_j$ holds for $i\ne j$ if and only if $(i,j)=(2,1)$ or $(3,1)$. In particular, the classes are pairwise distinct. The proof translates universal inclusion into a spanning-embedding question for coprimality graphs on consecutive integer intervals. Edge-count identities, arithmetic-progression estimates, and prime-residue switching exclude all other inclusions; a prime-gap argument reduces the remaining finite range to fifteen explicit rational inequalities. We further show that, for each nonexceptional ordered pair, the graph orders at which edge counts obstruct inclusion have a positive natural density. This density is described through a limiting distribution of correlated Euler products, for which we give an atomlessness argument. The classification resolves the distinct-index formulation of the universal-implication conjecture.
\end{abstract}
\section{Introduction}
A $k$-prime labeling of a graph $G$ of order $n$ is a bijection $f:V(G)\to\{k,k+1,\ldots,k+n-1\}$ such that $\gcd(f(u),f(v))=1$ for every edge $uv$. A graph admitting such a labeling is $k$-prime. This notion was introduced by Vaidya and Prajapati~\cite{VaidyaPrajapati}. Zhang and Li introduced the term \emph{uniform prime labeling} in their work on trees~\cite{ZhangLiTrees} and studied implications among $k$-prime labelings in their cycle-snake paper~\cite{ZhangLiSnakes}. Earlier investigations include small-order $k$-prime graphs~\cite{AbughneimAbughazaleh} and coprime graphs on consecutive integer intervals~\cite{BerkoveBrilleslyper,BaniMostafaGhorbani}; see also the graph-labeling survey~\cite{Gallian}. The latter coprimality graphs provide a natural model for our question, but their basic modular periodicity is already known. The following are the two conjectures formulated by Zhang and Li~\cite{ZhangLiSnakes}, numbered 1 and 2 in that work.
\begin{conjecture}\label{conj:distinct}
For distinct positive integers $i,j$, the classes of $i$-prime and $j$-prime graphs are different.
\end{conjecture}
\begin{conjecture}\label{conj:contain}
For distinct positive integers $i,j$, if every $i$-prime graph is $j$-prime, then $(i,j)\in\{(2,1),(3,1)\}$.
\end{conjecture}
Conjecture~\ref{conj:contain} is the corrected distinct-index formulation of the original Conjecture~2 in~\cite{ZhangLiSnakes}. Without $i\ne j$, the proposed implication classification would be false because $\mathcal P_i\subseteq\mathcal P_i$ for every $i$. Our classification establishes precisely this corrected formulation. Our main theorem resolves Conjecture~\ref{conj:contain} of~\cite{ZhangLiSnakes}; Conjecture~\ref{conj:distinct} follows as an immediate corollary. The proof proceeds by reducing universal implications to spanning embeddings of interval coprimality graphs. An edge-count translation identity and a conditional mean-value formula yield arithmetic criteria for nonimplication. Large prime divisors and weighted prime-residue switches cover all but finitely many ordered pairs; Nagura's prime-interval theorem bounds the remaining range, and an elementary prime-gap classification reduces the remaining cases to fifteen explicit rational inequalities. The additional finite classifications record first separating orders detected by the stated invariants; they are not prerequisites for the main theorem. The finite classification is proved using a short prime-gap table; the verification programs provide independent checks. The probabilistic section uses established methods of multiplicative-function distribution theory rather than claiming a new general limit theorem. We use standard facts on primes, totients, and divisor sums; see~\cite{HardyWright,Apostol} for background.
The principal graph-theoretic result is the complete inclusion classification, stated first in Section~\ref{sec:main-classification}. Its proof is given there in terms of later, fully proved arithmetic criteria. The coprimality-graph framework and nonimplication arguments are developed next; the quantitative density results follow their prerequisites. The finite certificates are displayed explicitly and the verification programs are included for reproducibility.
\section{Main classification}\label{sec:main-classification}
For each positive integer $k$, let $\mathcal P_k$ denote the class of all finite $k$-prime graphs.
\begin{lemma}[Known positive implications]\label{lem:positive}
Every $2$-prime graph and every $3$-prime graph is $1$-prime.
\end{lemma}
\begin{proof}
Let $G$ have $n$ vertices. If $n=1$, the assertion is immediate. Given a $2$-prime labeling, replace its label $n+1$ by $1$ and leave the other labels unchanged. This produces the label set $\{1,2,\ldots,n\}$; since $1$ is coprime to every integer, all edges remain valid.

Now suppose $n\ge2$ and let $f$ be a $3$-prime labeling. Replace the labels $n+1,n+2$ by $1,2$, respectively, when $n$ is even, and by $2,1$, respectively, when $n$ is odd. The label replaced by $2$ is even in either case. Thus each new label divides the old label at its vertex. If $uv$ is an edge, then $\gcd(f(u),f(v))=1$ implies that the new labels at $u,v$ are also coprime. The resulting labeling uses exactly $\{1,\ldots,n\}$ and is $1$-prime. These implications were also established in~\cite{ZhangLiSnakes}.
\end{proof}
\begin{theorem}[Complete classification of universal implications]\label{thm:classification}
For distinct positive integers $i,j$,
\[
 \mathcal P_i\subseteq\mathcal P_j
 \quad\Longleftrightarrow\quad
 (i,j)\in\{(2,1),(3,1)\}.
\]
\end{theorem}
\begin{proof}
The two stated implications follow from Lemma~\ref{lem:positive}.
Suppose $(i,j)$ is not one of these pairs. If $i<j$, apply Theorem~\ref{thm:increasing}. If $i>j$ and $j=1$, then $i\ge4$ because the two exceptional cases were excluded, so apply Theorem~\ref{thm:targetone}; if $j=2$ or $3$, apply Corollary~\ref{cor:targets23}. If $i>j\ge4$ and $i\le2j$, apply Corollary~\ref{cor:twice}. In the remaining range put $r=i-j>j\ge4$. Theorem~\ref{thm:large-r} handles $r\ge34$; Proposition~\ref{prop:finite435} handles $5\le r<34$ using the prime-gap characterization and the fifteen exact switching inequalities. These cases exhaust all possibilities.
\end{proof}
\begin{corollary}[Resolution of Conjecture 1 in~\cite{ZhangLiSnakes}]\label{cor:distinct}
For all distinct positive integers $i,j$, one has $\mathcal P_i\ne\mathcal P_j$.
\end{corollary}
\begin{proof}
Equality would imply both inclusions. By Theorem~\ref{thm:classification}, neither exceptional ordered pair admits the reverse inclusion.
\end{proof}

\section{Coprimality graphs and universal implications}
Let $H_{k,n}$ denote the graph with vertex set $\{k,\ldots,k+n-1\}$, with distinct vertices adjacent exactly when they are coprime. This is the $k$-coprime graph studied by Bani Mostafa and Ghorbani~\cite{BaniMostafaGhorbani}.
\begin{theorem}\label{thm:criterion}
For positive integers $i,j$, the following are equivalent:
\begin{enumerate}
\item Every $i$-prime graph is $j$-prime.
\item For every $n\geq1$, $H_{i,n}$ is $j$-prime.
\item For every $n\geq1$, there exists a bijection $\phi_n:\{i,\ldots,i+n-1\}\to\{j,\ldots,j+n-1\}$ carrying each coprime pair to a coprime pair.
\end{enumerate}
\end{theorem}
\begin{proof}
The identity labeling makes $H_{i,n}$ $i$-prime, so (1) implies (2). A $j$-prime labeling of $H_{i,n}$ is precisely a bijection from its label set to $\{j,\ldots,j+n-1\}$ that sends every edge to a coprime pair; hence (2) and (3) are equivalent. Finally, suppose (2) holds and $G$ has an $i$-prime labeling $f$. Identifying each vertex $v$ with $f(v)$ realizes $G$ as a spanning subgraph of $H_{i,n}$. If $g$ is a $j$-prime labeling of $H_{i,n}$, then $g\circ f$ is a bijection onto $\{j,\ldots,j+n-1\}$ and is coprime on every edge of $G$. Thus (2) implies (1).
\end{proof}
\begin{theorem}[Fixed-order periodicity]\label{thm:period}
For $n\geq2$, put $L_n=\operatorname{lcm}(1,\ldots,n-1)$. If $a\equiv b\pmod{L_n}$, then $H_{a,n}$ and $H_{b,n}$ are isomorphic under translation of label positions. In particular, a graph of order $n$ is $a$-prime if and only if it is $b$-prime. This period is generally not optimal; compare the primorial periodicity in~\cite{BerkoveBrilleslyper}.
\end{theorem}
\begin{proof}
For $0\leq r<s<n$, set $d=s-r$. Since $d\mid L_n$ and $a\equiv b\pmod d$, we have $\gcd(a+r,a+s)=\gcd(a+r,d)=\gcd(b+r,d)=\gcd(b+r,b+s)$. Thus adjacency is preserved.
\end{proof}
\begin{theorem}[Optimal squarefree translation period]\label{thm:primorial}
For $n\geq2$, put
\[
 P_n=\prod_{\substack{p<n\\p\text{ prime}}}p.
\]
If $a\equiv b\pmod{P_n}$, then $H_{a,n}$ and $H_{b,n}$ have identical adjacency relations under the position-preserving translation $a+r\mapsto b+r$. Moreover, $P_n$ is the least positive common translation period: if a positive integer $T$ has this property for every starting index $a$, then $P_n\mid T$. The existence of the primorial period is known~\cite{BerkoveBrilleslyper}; the least-period property also follows from their count of distinct modular tables. We include a direct proof for completeness.
\end{theorem}
\begin{proof}
For $0\leq r<s<n$, the pair of positions $r,s$ is adjacent precisely when no prime divisor of $s-r$ divides $a+r$. Every such prime is less than $n$, so congruence modulo $P_n$ preserves adjacency.
Conversely, let $T$ preserve all position-wise adjacencies for every starting index. Fix any prime $p<n$ and take $a=p$. Positions $0,p$ are nonadjacent because $\gcd(p,2p)=p$. Their translates must also be nonadjacent, so $\gcd(p+T,2p+T)=\gcd(p+T,p)>1$. Hence $p\mid T$. As this holds for every prime $p<n$, $P_n\mid T$.
\end{proof}
\begin{corollary}\label{cor:squarefree-lift}
If $H_{i,n}$ is not $j$-prime and $a\equiv i\pmod{P_n}$, $b\equiv j\pmod{P_n}$, then $a$-primality does not universally imply $b$-primality.
\end{corollary}
\begin{proof}
Apply Theorem~\ref{thm:primorial} to the source and target graphs, followed by Theorem~\ref{thm:criterion}.
\end{proof}
\begin{corollary}\label{cor:lift}
If $H_{i,n}$ is not $j$-prime and $a\equiv i\pmod{L_n}$, $b\equiv j\pmod{L_n}$, then $a$-primality does not universally imply $b$-primality.
\end{corollary}
\begin{proof}
Theorem~\ref{thm:period} identifies the respective coprimality graphs. Apply Theorem~\ref{thm:criterion}.
\end{proof}
\section{Distinctness of all $k$-prime graph classes}
For positive integers $x$ and $m$, put
\[
 C_m(x)=\#\{1\le d\le m:\gcd(d,x)=1\},\qquad e(k,n)=|E(H_{k,n})|.
\]
\begin{lemma}[Translation identity]\label{lem:translation}
For $k\ge1$ and $n\ge2$,
\[
 e(k+1,n)-e(k,n)=C_{n-1}(n+k)-C_{n-1}(k).
\]
Consequently, for $i<j$,
\[
 e(j,n)-e(i,n)=\sum_{k=i}^{j-1}\bigl(C_{n-1}(n+k)-C_{n-1}(k)\bigr).
\]
\end{lemma}
\begin{proof}
The intervals $[k,k+n-1]$ and $[k+1,k+n]$ have $n-1$ labels in common. Edges between common labels cancel. The removed label $k$ has $C_{n-1}(k)$ neighbors and the added label $k+n$ has $C_{n-1}(k+n)$ neighbors. Summing the resulting identity gives the second assertion.
\end{proof}
\begin{lemma}\label{lem:coprime-count}
For all integers $q,m\ge1$,
\[
 C_m(q)=m\frac{\varphi(q)}q+O(\tau(q)),
\]
with an absolute implied constant, where $\tau(q)$ denotes the number of positive divisors of $q$. In particular, the error is $O(q)$. If $Q\mid x$, then $C_m(x)\le C_m(Q)$.
\end{lemma}
\begin{proof}
By M\"obius inversion,
\[
 C_m(q)=\sum_{d\mid q}\mu(d)\left\lfloor\frac md\right\rfloor
 =m\sum_{d\mid q}\frac{\mu(d)}d+O\!\left(\sum_{d\mid q}1\right)
 =m\frac{\varphi(q)}q+O(\tau(q)).
\]
The bound $\tau(q)\le q$ gives the weaker estimate. If $Q\mid x$, every integer coprime to $x$ is also coprime to $Q$, giving the final assertion.
\end{proof}
\begin{theorem}[Nonimplication in the increasing-index direction]\label{thm:increasing}
For every $1\le i<j$, there are a positive integer $Q$, a residue class $a\pmod Q$, and a constant $c>0$ such that
\[
 e(i,n)-e(j,n)\ge cn
\]
for every sufficiently large positive $n\equiv a\pmod Q$. In particular, $i$-primality does not imply $j$-primality.
\end{theorem}
\begin{proof}
Put $r=j-i$ and $A=\sum_{k=i}^{j-1}\varphi(k)/k>0$. Choose $\epsilon>0$ such that $r\epsilon<A$. For each $k=i,\ldots,j-1$, choose a squarefree positive integer $Q_k$ with $\varphi(Q_k)/Q_k<\epsilon$, arranging that the $Q_k$ are pairwise coprime. This is possible by partitioning sufficiently many primes into $r$ disjoint finite sets whose respective products of $(1-1/p)$ are less than $\epsilon$; it uses divergence of the sum of reciprocals of primes~\cite{HardyWright}.
By the Chinese remainder theorem there is a residue class $a\pmod Q$, $Q=\prod_k Q_k$, such that $a\equiv-k\pmod{Q_k}$ for each $k$. For $n\equiv a\pmod Q$, each $Q_k$ divides $n+k$, so Lemma~\ref{lem:coprime-count} yields
\[
 \sum_{k=i}^{j-1}C_{n-1}(n+k)
 \le r\epsilon(n-1)+O_{i,j,Q}(1).
\]
The same lemma gives
\[
 \sum_{k=i}^{j-1}C_{n-1}(k)=A(n-1)+O_{i,j}(1).
\]
By Lemma~\ref{lem:translation},
\[
 e(i,n)-e(j,n)\ge (A-r\epsilon)(n-1)-O_{i,j,Q}(1),
\]
which is at least $cn$ for all sufficiently large admissible $n$, for a suitable $c>0$. The edge-count obstruction finishes the proof.
\end{proof}

\section{Arithmetic and graph-theoretic obstructions}
\begin{proposition}\label{prop:odd-even}
If $i$ is odd and $j$ is even, then $i$-primality does not imply $j$-primality.
\end{proposition}
\begin{proof}
The triangle $K_3$ is $i$-prime, since $i,i+1,i+2$ are pairwise coprime when $i$ is odd. It is not $j$-prime, since $j$ and $j+2$ are both even.
\end{proof}
\begin{theorem}[Prime-divisor separation]\label{thm:prime}
Let $p$ be prime and $L=\operatorname{lcm}(1,\ldots,p-1)$. If $i\equiv j\pmod L$, $p\nmid i$, and $p\mid j$, then $i$-primality does not imply $j$-primality.
\end{theorem}
\begin{proof}
In $H_{i,p+1}$ and $H_{j,p+1}$, corresponding pairs at distance $d<p$ have the same coprimality status because $d\mid L$. The unique pair at distance $p$ is coprime in the first graph and noncoprime in the second. Hence $|E(H_{i,p+1})|=|E(H_{j,p+1})|+1$, which prevents a spanning embedding of the former into the latter.
\end{proof}
\begin{corollary}\label{cor:infinitely-many}
For each fixed positive integer $i$, infinitely many positive integers $j$ satisfy that $i$-primality does not imply $j$-primality.
\end{corollary}
\begin{proof}
Choose a prime $p\nmid i$. Since $\gcd(p,L)=1$, the Chinese remainder theorem supplies infinitely many positive $j$ with $j\equiv i\pmod L$ and $j\equiv0\pmod p$. Apply Theorem~\ref{thm:prime}.
\end{proof}
\begin{proposition}[Embedding obstructions]\label{prop:invariants}
Let $e_{k,n}=|E(H_{k,n})|$, and write $d_1^{(k)}\leq\cdots\leq d_n^{(k)}$ for the sorted degrees of $H_{k,n}$. If for some $n$ either $e_{i,n}>e_{j,n}$ or $d_r^{(i)}>d_r^{(j)}$ for some $r$, then $i$-primality does not imply $j$-primality.
\end{proposition}
\begin{proof}
Any spanning embedding of $H_{i,n}$ into $H_{j,n}$ preserves all source edges, hence cannot increase the source edge count beyond the target edge count. It also gives a bijection matching each source vertex to a target vertex of at least its degree. Such a matching requires coordinatewise inequality of the sorted degree sequences.
\end{proof}
\section{Nonimplications in every arithmetic progression}
The prime-divisor obstruction yields a uniform conclusion that does not require fixing the parity or a particular congruence class of the indices.
\begin{theorem}[Arithmetic-progression separation]\label{thm:progression}
For every pair of positive integers $i,M$, there are infinitely many positive integers $j$ with
\[
 j\equiv i\pmod M
\]
such that $i$-primality does not universally imply $j$-primality. More precisely, for some positive integer $Q$, all positive integers in one residue class modulo $Q$ have these properties.
\end{theorem}
\begin{proof}
Choose a prime $p>M$ with $p\nmid i$, and put $L=\operatorname{lcm}(1,\ldots,p-1)$. Because $M<p$, we have $M\mid L$, while $\gcd(p,L)=1$. The Chinese remainder theorem supplies a residue class $j_0\pmod{pL}$ satisfying
\[
 j_0\equiv i\pmod L,\qquad j_0\equiv0\pmod p.
\]
Every positive integer $j\equiv j_0\pmod{pL}$ satisfies $j\equiv i\pmod M$ and, by Theorem~\ref{thm:prime}, gives the claimed nonimplication. Take $Q=pL$.
\end{proof}
\begin{corollary}\label{cor:density}
For fixed positive integers $i,M$, the set of positive integers $j\equiv i\pmod M$ for which $i$-primality does not imply $j$-primality has positive lower natural density among the positive integers. Specifically, for any choice of $p$ as in the proof of Theorem~\ref{thm:progression}, this lower density is at least $1/(pL)$, where $L=\operatorname{lcm}(1,\ldots,p-1)$.
\end{corollary}
\begin{proof}
The set contains the complete residue class $j_0\pmod{pL}$ constructed above, which has natural density $1/(pL)$.
\end{proof}
\begin{corollary}\label{cor:sameparity}
For every positive integer $i$, infinitely many $j\equiv i\pmod6$ satisfy $i$-prime $\not\Rightarrow$ $j$-prime universally.
\end{corollary}
\begin{proof}
Apply Theorem~\ref{thm:progression} with $M=6$.
\end{proof}

\section{Arbitrarily delayed separating graphs}
\begin{theorem}[Exact first separating order]\label{thm:delayed}
Let $i$ be a positive integer and let $p$ be a prime with $p\nmid i$.
There exist infinitely many positive integers $j$ for which every graph
of order at most $p$ is $i$-prime if and only if it is $j$-prime, but
$H_{i,p+1}$ is $i$-prime and not $j$-prime. In particular, the smallest
order of a graph that is $i$-prime but not $j$-prime is exactly $p+1$.
\end{theorem}
\begin{proof}
Put $L=\operatorname{lcm}(1,\ldots,p-1)$. Since $\gcd(p,L)=1$,
the Chinese remainder theorem supplies infinitely many positive integers
$j$ satisfying $j\equiv i\pmod L$ and $j\equiv0\pmod p$.
For $2\le n\le p$, the period $L_n$ divides $L$, so
Theorem~\ref{thm:period} identifies $H_{i,n}$ and $H_{j,n}$ under
translation of positions. The same assertion for $n=1$ is immediate.
Therefore the two primality notions agree on all graphs of order at most
$p$. At order $p+1$, corresponding pairs of positions at distance less
than $p$ have the same coprimality status. The unique pair at distance
$p$ is coprime for the initial index $i$ and noncoprime for the initial
index $j$. Hence $H_{i,p+1}$ has exactly one more edge than $H_{j,p+1}$,
so it cannot admit a $j$-prime labeling.
\end{proof}
\begin{corollary}\label{cor:no-uniform-bound}
For every positive integer $i$ and every integer $N\ge2$, there are
infinitely many distinct indices $j$ such that $i$-primality and
$j$-primality agree on every graph of order at most $N$, but
$i$-primality does not universally imply $j$-primality.
\end{corollary}
\begin{proof}
Choose a prime $p\ge N$ not dividing $i$ and apply
Theorem~\ref{thm:delayed}.
\end{proof}

\section{Conditional coprimality densities and reverse nonimplications}
For squarefree $M\ge1$ and any integer $b$, define
\[
 D_M(b)=\prod_{p\nmid M}\left(1-\frac1{p^2}\right)
 \prod_{\substack{p\mid M\\p\mid b}}\left(1-\frac1p\right).
\]
\begin{lemma}[Arithmetic-progression mean]\label{lem:conditional}
For fixed squarefree $M$, residue $a\pmod M$, and positive integer $k$, as $X\to\infty$,
\[
 \sum_{\substack{X\le n\le2X\\n\equiv a\pmod M}}C_{n-1}(n+k)
 =D_M(a+k)\sum_{\substack{X\le n\le2X\\n\equiv a\pmod M}}(n-1)
 +O_{M,k}(X\log(2X)).
\]
\end{lemma}
\begin{proof}
M\"obius inversion gives
\[
 C_{n-1}(n+k)=(n-1)\sum_{d\mid n+k}\frac{\mu(d)}d+O(\tau(n+k)).
\]
The classical divisor-sum estimate~\cite{Apostol} shows that the sum of the divisor-function errors for $X\le n\le2X$ is $O_k(X\log(2X))$. For squarefree $d$ let $g=\gcd(d,M)$. The congruences $n\equiv a\pmod M$ and $n\equiv-k\pmod d$ are compatible exactly when $g\mid a+k$. In that case the common solutions form one residue class modulo $L=\operatorname{lcm}(d,M)=Md/g$. For any such residue class, the elementary progression estimate gives
\[
 \sum_{\substack{X\le n\le2X\\n\equiv c\pmod L}}(n-1)
 =\frac{3X^2}{2L}+O(X+1).
\]
Because $1/L=g/(Md)$, comparison with the progression modulo $M$ shows that the weighted count for the simultaneous congruences is $g/d$ times the weighted count for $n\equiv a\pmod M$, with error $O(X+1)$ for each $d$. Multiplication by $|\mu(d)|/d$ and summation over $d\le 2X+k$ give a total error bounded by
\[
 O\!\left((X+1)\sum_{d\le2X+k}\frac1d\right)
 =O_k(X\log(2X)).
\]
The implied constant may depend on the fixed modulus $M$ and shift $k$. The absolutely convergent main coefficient is
\[
 \sum_{\substack{d\ge1\\\gcd(d,M)\mid a+k}}
 \frac{\mu(d)\gcd(d,M)}{d^2}=D_M(a+k).
\]
The discarded tail contributes $O_{M,k}(X)$, since $\gcd(d,M)\le M$. This proves the formula.
\end{proof}
\begin{theorem}[Conditional-density obstruction]\label{thm:conditional}
Let $1\le j<i$ and suppose that for some squarefree $M$ and residue $a\pmod M$,
\[
 \sum_{k=j}^{i-1}D_M(a+k)>
 \sum_{k=j}^{i-1}\frac{\varphi(k)}k.
\]
Then $e(i,n)>e(j,n)$ for infinitely many $n\equiv a\pmod M$, and hence $i$-primality does not imply $j$-primality.
\end{theorem}
\begin{proof}
Sum Lemma~\ref{lem:translation} over $n\in[X,2X]$ in the given residue class. Apply Lemma~\ref{lem:conditional} to each shifted term, and Lemma~\ref{lem:coprime-count} to the fixed terms. The resulting sum of $e(i,n)-e(j,n)$ equals a strictly positive constant times $\sum(n-1)$ plus $O_{i,j,M}(X\log(2X))$. Writing the strict difference of the two coefficient sums as $\delta>0$, the total is
\[
 \delta\sum_{\substack{X\le n\le2X\\n\equiv a\pmod M}}(n-1)
 +O_{i,j,M}(X\log(2X))
 =\frac{3\delta}{2M}X^2+O_{i,j,M}(X\log(2X)).
\]
The main term is positive and of order $X^2$, whereas the error is $o(X^2)$. Therefore some $n$ in each sufficiently large dyadic interval has $e(i,n)>e(j,n)$. Choosing disjoint dyadic intervals gives infinitely many distinct such orders. By Theorem~\ref{thm:criterion}, the graph $H_{i,n}$ then witnesses failure of the universal implication.
\end{proof}

\section{Densities of separating graph orders}\label{sec:order-density}
\begin{theorem}[Positive-density separation orders]\label{thm:order-density}
For distinct positive integers $i,j$, put
\[
 S_{i,j}=\{n\ge2:e(i,n)>e(j,n)\}.
\]
Then $S_{i,j}$ has positive lower natural density if and only if
$(i,j)\notin\{(2,1),(3,1)\}$. More precisely, if $i<j$,
$S_{i,j}$ contains all sufficiently large members of a fixed arithmetic
progression. If $i>j$ and $(i,j)$ is not exceptional, there are a squarefree
integer $M$, a residue $a\pmod M$, and a number $\delta>0$ such that
\[
 \underline d(S_{i,j})\ge \frac{\delta}{2M(i-j)}>0,
\]\[
 \delta=\sum_{k=j}^{i-1}D_M(a+k)
       -\sum_{k=j}^{i-1}\frac{\varphi(k)}k.
\]
\end{theorem}
\begin{proof}
For $i<j$, Theorem~\ref{thm:increasing} supplies a fixed arithmetic
progression on which $e(i,n)>e(j,n)$ for all sufficiently large $n$.
Its density is positive.

Suppose $i>j$ and $(i,j)$ is not exceptional. The proofs of
Theorem~\ref{thm:targetone}, Corollary~\ref{cor:targets23},
Corollary~\ref{cor:twice}, Theorem~\ref{thm:large-r}, and
Proposition~\ref{prop:finite435} each establish the strict hypothesis
of Theorem~\ref{thm:conditional} for some fixed squarefree $M$ and
residue $a\pmod M$. Fix such a choice and let $r=i-j$ and $\delta>0$
be as in the statement.

By Lemmas~\ref{lem:translation} and~\ref{lem:conditional}, together with
Lemma~\ref{lem:coprime-count}, applying the arithmetic-progression mean estimate on successive dyadic
intervals, summing their errors, and absorbing the finitely many initial
terms gives
\[
\begin{aligned}
 &\sum_{\substack{2\le n\le X\\n\equiv a\pmod M}}
     \bigl(e(i,n)-e(j,n)\bigr)\\
 &\quad=\delta\sum_{\substack{2\le n\le X\\n\equiv a\pmod M}}(n-1)
     +O_{i,j,M}(X\log^2(2X))\\
 &\quad=\frac{\delta}{2M}X^2+O_{i,j,M}(X\log^2(2X)).
\end{aligned}
\]
The translation identity also gives
$e(i,n)-e(j,n)\le r(n-1)\le rX$ for $n\le X$.
Discarding nonpositive summands therefore yields
\[
 rX\,\#\{n\le X:e(i,n)>e(j,n)\}
 \ge \frac{\delta}{2M}X^2+O_{i,j,M}(X\log^2(2X)).
\]
Division by $rX^2$ and passage to the lower limit prove the bound.
Finally, for $(i,j)=(2,1)$ or $(3,1)$, Lemma~\ref{lem:positive}
forces $e(i,n)\le e(j,n)$ for every $n$, since $H_{i,n}$ admits
an embedding into $H_{j,n}$. Thus $S_{i,j}$ is empty.
\end{proof}

\begin{corollary}[Relative density in a separating progression]\label{cor:relative-order-density}
Let $i>j$ be distinct positive integers with
$(i,j)\notin\{(2,1),(3,1)\}$. Choose a squarefree integer $M$ and a
residue $a\pmod M$ for which
\[
\delta=\sum_{k=j}^{i-1}D_M(a+k)
-\sum_{k=j}^{i-1}\frac{\varphi(k)}k>0.
\]
Then
\[
\liminf_{X\to\infty}
\frac{\#\{2\le n\le X:n\equiv a\pmod M,\ e(i,n)>e(j,n)\}}
{\#\{2\le n\le X:n\equiv a\pmod M\}}
\ge\frac{\delta}{2(i-j)}>0.
\]
\end{corollary}
\begin{proof}
Put $r=i-j$. The proof of Theorem~\ref{thm:order-density} gives
\[
\sum_{\substack{2\le n\le X\\n\equiv a\pmod M}}
(e(i,n)-e(j,n))
=\frac{\delta}{2M}X^2+O_{i,j,M}(X\log^2(2X)).
\]
Every positive summand is at most $rX$, while nonpositive summands
can only decrease the sum. Hence the number of separating orders in
this progression is at least
$\delta X/(2Mr)+O_{i,j,M}(\log^2(2X))$.
Since the progression contains $X/M+O_M(1)$ integers up to $X$,
dividing and taking the lower limit gives the result.
\end{proof}

\section{Limiting distributions and natural densities of separating orders}\label{sec:natural-density}
We strengthen the lower-density result of Theorem~\ref{thm:order-density} to the existence of a natural density. The probabilistic argument below is independent of the finite certificates once the positive lower-density theorem is established.
The existence and continuity of the limiting distribution of the individual totient ratio $\varphi(n)/n$ go back to Schoenberg~\cite{Schoenberg1928}; see also Erd\H{o}s~\cite{Erdos1974} and Weingartner~\cite{Weingartner2007} for related distributional work. The truncation of multiplicative functions to small primes and the construction of limiting laws from independent prime components are also classical; see the general infinite-convolution framework of Jessen and Wintner~\cite{JessenWintner1935}. Our specific application concerns the correlated sum of several consecutive shifts and its sign distribution. We do not claim the underlying probabilistic or concentration techniques as new; the argument here concerns the specific correlated Euler-product sum arising from consecutive shifts.

Fix integers $i>j\ge1$, put $r=i-j$, and set $K=\{j,j+1,\ldots,i-1\}$. Write
\[
 A_{j,i}=\sum_{k\in K}\frac{\varphi(k)}k,\qquad
 F_{i,j}(n)=\sum_{k\in K}\frac{\varphi(n+k)}{n+k}-A_{j,i}.
\]
For each prime $p$, let $U_p$ be an independent uniformly distributed residue modulo $p$. Define
\[
 Z_k=\prod_p\left(1-\frac1p\right)^{\mathbf1_{\{U_p\equiv-k\pmod p\}}},\qquad
 Z_{i,j}=\sum_{k\in K}Z_k-A_{j,i}.
\]
All these products converge almost surely to strictly positive random variables: for each $k$,
\[
 \sum_p\mathbb E\left[\frac{\mathbf1_{\{U_p\equiv-k\pmod p\}}}{p}\right]
 =\sum_p\frac1{p^2}<\infty,
\]
and $-\log(1-1/p)=O(1/p)$ for $p\ge3$; the factor at $p=2$ is harmless.

\begin{theorem}[Limiting Euler-product distribution]\label{thm:limit-distribution}
For fixed $i>j\ge1$, the empirical distributions of $F_{i,j}(n)$ converge weakly to the distribution of $Z_{i,j}$. Moreover,
\[
 \frac{e(i,n)-e(j,n)}n=F_{i,j}(n)+O_{i,j}(n^{-1/2})
 \qquad(n\to\infty).
\]
\end{theorem}
\begin{proof}
By the edge-count translation identity (Lemma~\ref{lem:translation}) and the elementary coprimality-count estimate (Lemma~\ref{lem:coprime-count}),
\[
 C_{n-1}(n+k)=(n-1)\frac{\varphi(n+k)}{n+k}+O(\tau(n+k)),
 \qquad C_{n-1}(k)=(n-1)\frac{\varphi(k)}k+O_k(1).
\]
For fixed $k$, $\tau(n+k)=O_k(\sqrt n)$, proving the displayed approximation.

For $y\ge2$, let $Q_y=\prod_{p\le y}p$ and
\[
 f_y(m)=\prod_{\substack{p\le y\\p\mid m}}(1-p^{-1}),\qquad
 F^{(y)}_{i,j}(n)=\sum_{k\in K}f_y(n+k)-A_{j,i}.
\]
The latter function is periodic modulo $Q_y$. Its empirical distribution therefore converges to the distribution of
\[
 Z^{(y)}_{i,j}=\sum_{k\in K}\prod_{p\le y}
 (1-p^{-1})^{\mathbf1_{\{U_p\equiv-k\pmod p\}}}-A_{j,i},
\]
by the Chinese remainder theorem. The elementary product inequality gives
\[
 0\le f_y(m)-\frac{\varphi(m)}m
 \le\sum_{\substack{p>y\\p\mid m}}\frac1p.
\]
For each fixed $k$, counting multiples of $p$ among $n+k$, $1\le n\le X$, yields
\[
 \limsup_{X\to\infty}\frac1X\sum_{n\le X}
 \left|f_y(n+k)-\frac{\varphi(n+k)}{n+k}\right|
 \le\sum_{p>y}\frac1{p^2}.
\]
Indeed, the finite-$X$ bound is at most $\sum_{y<p\le X+k}(p^{-2}+1/(Xp))$. The extra term is at most $X^{-1}\sum_{m\le X+k}m^{-1}=O_k(\log(2X)/X)$, which tends to zero. Hence the corresponding mean error for $F^{(y)}_{i,j}$ is at most $r\sum_{p>y}p^{-2}$ in the limit. On the random side, the truncated products converge almost surely to the positive products $Z_k$, and are bounded by one; therefore $Z^{(y)}_{i,j}\to Z_{i,j}$ in $L^1$. Applying these estimates to bounded Lipschitz test functions and then letting $y\to\infty$ proves weak convergence.
\end{proof}

\begin{lemma}[Distinct finite Euler products]\label{lem:subset-products}
If $F$ is a finite set of distinct primes, the numbers
\[
 \prod_{p\in A}\frac{p-1}{p},\qquad A\subseteq F,
\]
are pairwise distinct.
\end{lemma}
\begin{proof}
Suppose two such products are equal. First cancel the identical factors corresponding to primes in the intersection of the two indexing subsets. Let $P$ be the largest prime in their symmetric difference, and suppose that its factor $(P-1)/P$ belongs to the first remaining product. Every remaining prime $q$ is at most $P$, and $P$ occurs as a denominator prime in only the first product; numerators $q-1$ for $q<P$ are also not divisible by $P$. The factor $(P-1)/P$ has $P$-adic valuation $-1$, while all factors corresponding to primes $q<P$ have $P$-adic valuation $0$. Hence the two rational products have different $P$-adic valuations, a contradiction.
\end{proof}

\begin{theorem}[Atomlessness of the limiting distribution]\label{thm:atomless}
For each fixed $i>j\ge1$, the distribution of $Z_{i,j}$ has no atoms. In particular, $\Pr(Z_{i,j}=0)=0$.
\end{theorem}
\begin{proof}
Choose a distinguished $k_0\in K$. Reveal the residues $U_p$ for $p\le r$. For each prime $p>r$, reveal whether $U_p$ equals $-k\pmod p$ for some $k\in K\setminus\{k_0\}$ and, if so, which $k$. Let $\mathcal F$ be the sigma-field generated by this information and let $\mathcal Q$ be the random set of primes $p>r$ for which no such hit occurs. Since the $r$ residues $-k\pmod p$ are distinct, conditional on $\mathcal F$ the variables
\[
 B_p=\mathbf1_{\{U_p\equiv-k_0\pmod p\}},\qquad p\in\mathcal Q,
\]
are independent Bernoulli variables with parameters $\theta_p=(p-r+1)^{-1}$. All the other $Z_k$, $k\ne k_0$, are $\mathcal F$-measurable. Moreover,
\[
 Z_{k_0}=C\prod_{p\in\mathcal Q}(1-p^{-1})^{B_p}
\]
for an $\mathcal F$-measurable random constant $C>0$ (including the revealed small-prime factors).

For each $p>r$, the probability that $p\notin\mathcal Q$ is $(r-1)/p$. Thus
\[
 \mathbb E\sum_{\substack{p>r\\p\notin\mathcal Q}}\frac1p
 =(r-1)\sum_{p>r}\frac1{p^2}<\infty.
\]
The nonnegative random sum on the left is therefore finite almost surely. By divergence of the sum of reciprocal primes,
\[
 \sum_{p\in\mathcal Q}\theta_p
 \ge\sum_{p\in\mathcal Q}\frac1p=\infty
 \qquad\text{almost surely}.
\]
Fix such a realization of $\mathcal F$. The remaining product
$Y=\prod_{p\in\mathcal Q}(1-p^{-1})^{B_p}$ is strictly positive almost surely, since
\[
 \sum_{p\in\mathcal Q}\mathbb E[B_p/p\mid\mathcal F]
 =\sum_{p\in\mathcal Q}\frac1{p(p-r+1)}<\infty.
\]
For a finite $E\subset\mathcal Q$, put $Y_E=\prod_{p\in E}(1-p^{-1})^{B_p}$. By Lemma~\ref{lem:subset-products}, distinct Bernoulli configurations give distinct values, so
\[
 \sup_x\Pr(Y_E=x\mid\mathcal F)
 =\prod_{p\in E}\max\{\theta_p,1-\theta_p\}.
\]
As $E$ increases through the eligible primes, the right side tends to zero, since $\theta_p\to0$ and $\sum_{p\in\mathcal Q}\theta_p=\infty$. Writing $Y=Y_ET_E$, the tail $T_E$ is strictly positive and conditionally independent of $Y_E$. For every real $x$,
\[
 \Pr(Y=x\mid\mathcal F)
 \le\sup_t\Pr(Y_E=t\mid\mathcal F).
\]
It follows that $Y$, and hence $CY$ plus the $\mathcal F$-measurable sum of the other $Z_k$, has no conditional atoms. Integrating the conditional probabilities gives $\Pr(Z_{i,j}=x)=0$ for every $x$.
\end{proof}

\begin{theorem}[Natural density of separating graph orders]\label{thm:natural-order-density}
For every pair of distinct positive integers $i,j$, the natural density
\[
 d(S_{i,j})=\lim_{X\to\infty}\frac1X
 \#\{2\le n\le X:e(i,n)>e(j,n)\}
\]
exists. It is strictly positive if and only if $(i,j)\notin\{(2,1),(3,1)\}$. For $i>j$ it equals $\Pr(Z_{i,j}>0)$, with $Z_{i,j}$ as defined above; for $i<j$ it equals $\Pr(Z_{j,i}<0)$.
\end{theorem}
\begin{proof}
For $i>j$, Theorem~\ref{thm:limit-distribution} gives weak convergence of the empirical distributions of $(e(i,n)-e(j,n))/n$ to the law of $Z_{i,j}$. By Theorem~\ref{thm:atomless}, the boundary of $(0,\infty)$ has probability zero, so the portmanteau theorem implies that the frequency of positive values converges to $\Pr(Z_{i,j}>0)$. For $i<j$, reverse the sign of the difference and apply the same argument to $(j,i)$. Finally, Theorem~\ref{thm:order-density} shows that the resulting density is positive for every nonexceptional pair; for the two exceptional pairs it is zero by Lemma~\ref{lem:positive}.
\end{proof}

\begin{proposition}[A finite prime-pattern certificate]\label{prop:primepattern}
Let $1\le j<i$, let $T$ be a positive integer, and set
$M_T=\prod_{p\le T}p$ and $R_T=\prod_{p>T}(1-p^{-2})$.
For a residue $a\pmod {M_T}$ put
\[
 W_T(a;j,i)=\sum_{k=j}^{i-1}
 \prod_{\substack{p\le T\\p\mid a+k}}(1-p^{-1}).
\]
If
\[
 R_T W_T(a;j,i)>\sum_{k=j}^{i-1}\frac{\varphi(k)}k,
\]
then $i$-primality does not imply $j$-primality.
\end{proposition}
\begin{proof}
The defining Euler product gives
$\sum_{k=j}^{i-1}D_{M_T}(a+k)=R_TW_T(a;j,i)$.
The conclusion follows from Theorem~\ref{thm:conditional}.
\end{proof}
\begin{proposition}\label{prop:11to4}
There exist infinitely many positive integers $n$ such that
$e(11,n)>e(4,n)$. In particular, $11$-primality does not
imply $4$-primality.
\end{proposition}
\begin{proof}
Choose $T=11$, $M_T=2310$, and $a=33$. For
$k=4,\ldots,10$, the local factors contributed by the primes
$2,3,5,7,11$ are, in order,
\[
 1,\quad\frac12,\quad\frac23,\quad\frac25,
 \quad1,\quad\frac27,\quad1.
\]
Their sum is $W_T(33;4,11)=1019/210$, whereas
\[
 \sum_{k=4}^{10}\frac{\varphi(k)}k=\frac{142}{35}.
\]
Furthermore,
\[
 R_T=\prod_{p>11}(1-p^{-2})
 \ge 1-\sum_{p>11}p^{-2}
 >1-\sum_{m=12}^{\infty}m^{-2}
 >\frac{10}{11}.
\]
Thus
\[
 R_TW_T(33;4,11)>\frac{1019}{231}
 >\frac{142}{35}.
\]
Proposition~\ref{prop:primepattern} and
Theorem~\ref{thm:conditional} establish the assertion.
\end{proof}

\begin{theorem}[Large-prime-factor obstruction]\label{thm:largeprime}
Let $1\le j<i$, $r=i-j$. If at least one integer in $\{j,\ldots,i-1\}$ has a prime divisor greater than $r$, then $i$-primality does not imply $j$-primality.
\end{theorem}
\begin{proof}
For $T>r$ let $M_T=\prod_{p\le T}p$. For each prime $p\le r$ prescribe $a_T\equiv0\pmod p$. For $r<p\le T$, choose $a_T$ to avoid the at most $r$ residues $-j,\ldots,-(i-1)$ modulo $p$. The Chinese remainder theorem combines these choices. Then
\[
 \sum_{k=j}^{i-1}D_{M_T}(a_T+k)
 =\prod_{p>T}(1-p^{-2})
 \sum_{k=j}^{i-1}\prod_{\substack{p\le r\\p\mid k}}(1-p^{-1}).
\]
As $T\to\infty$, the prefactor tends to $1$. Every summand in the remaining sum is at least $\varphi(k)/k$, with a strict inequality for an integer having a prime divisor greater than $r$. Hence for sufficiently large finite $T$, Theorem~\ref{thm:conditional} applies.
\end{proof}
\begin{corollary}\label{cor:twice}
For all $2\le j<i\le2j$, $i$-primality does not imply $j$-primality.
\end{corollary}
\begin{proof}
Put $r=i-j\le j$. If $j>r$, the Sylvester--Schur theorem~\cite{Erdos1934} gives a prime divisor greater than $r$ in the block $j,\ldots,j+r-1$. If $j=r$, Bertrand's postulate~\cite{HardyWright} supplies a prime $p$ with $j<p<2j$, which lies in the block. Apply Theorem~\ref{thm:largeprime}.
\end{proof}
\begin{lemma}[Shifted totient-sum criterion]\label{lem:shifted}
Let $1\le j<i$ and $r=i-j$. If
\[
 \sum_{m=1}^{r}\frac{\varphi(m)}m>
 \sum_{m=j}^{i-1}\frac{\varphi(m)}m,
\]
then $i$-primality does not imply $j$-primality.
\end{lemma}
\begin{proof}
For primes $p\le r$ impose $a+j\equiv1\pmod p$. For $r<p\le T$, avoid all $r$ residues $-j,\ldots,-(i-1)$ modulo $p$. The Chinese remainder theorem yields a residue modulo $M_T=\prod_{p\le T}p$. As $T\to\infty$, the sum of conditional densities tends to $\sum_{m=1}^r\varphi(m)/m$. The strict hypothesis and Theorem~\ref{thm:conditional} finish the proof.
\end{proof}
\begin{corollary}\label{cor:targets23}
For every $i\ge3$, $i$-primality does not imply $2$-primality. For every $i\ge4$, $i$-primality does not imply $3$-primality.
\end{corollary}
\begin{proof}
Write $w(m)=\varphi(m)/m$. For $j=2$ the difference between the two sums in Lemma~\ref{lem:shifted} is $w(1)-w(r+1)>0$. For $j=3$ it is $w(1)+w(2)-w(r+1)-w(r+2)>0$, since the last two integers are consecutive and at least $2$, one is even, and both have totient ratio strictly below $1$.
\end{proof}

\begin{theorem}[Nonimplication to ordinary prime labeling]\label{thm:targetone}
For every integer $i\ge4$, there are infinitely many positive integers $n$ for which
\[
 e(i,n)>e(1,n).
\]
In particular, $i$-primality does not imply $1$-primality.
\end{theorem}
\begin{proof}
Put $r=i-1\ge3$ and $w(m)=\varphi(m)/m$. By Bertrand's
postulate there is an odd prime $p$ with $(r+1)/2<p\le r$.
Indeed, for odd $r=2h+1$, apply Bertrand's postulate to $h+1$;
for even $r=2h$, apply it to $h$, noting that $r\ge4$.
Thus $p-1\in\{1,\ldots,r\}$ and $2p>r+1$.

For $T>r$, let $M_T=\prod_{q\le T}q$, the product over primes.
Choose $a_T\pmod{M_T}$ by the Chinese remainder theorem so that
\[
\begin{cases}
 a_T\equiv0\pmod q,&q\le r,\ q\ne p,\\
 a_T\equiv1\pmod p,&q=p,\\
 a_T\not\equiv-m\pmod q\quad(1\le m\le r),&r<q\le T.
\end{cases}
\]
The final conditions are feasible because each $q>r$ has more
residue classes than the $r$ forbidden values.
Let
\[
 v(m)=\prod_{\substack{q\le r,\ q\ne p\\q\mid m}}
 (1-q^{-1}).
\]
For the original block $1,\ldots,r$, the unique multiple of $p$
is $p$, and $v(p)=1$. For the shifted block $a_T+1,\ldots,a_T+r$,
the unique multiple of $p$ is at position $p-1$.
Consequently, writing $R_T=\prod_{q>T}(1-q^{-2})$, the conditional
density formula yields
\[
\sum_{m=1}^r D_{M_T}(a_T+m)
 =R_T\left(\sum_{m=1}^r w(m)
 +\frac{1-v(p-1)}p\right).
\]
Since $p-1$ is even, $v(p-1)\le1/2$, so the expression in
parentheses exceeds $\sum_{m=1}^r w(m)$ by at least $1/(2p)$.
As $T\to\infty$, $R_T\to1$. Thus, for sufficiently large finite
$T$, the strict hypothesis of Theorem~\ref{thm:conditional}
holds with $j=1$ and $i=r+1$. The result follows.
\end{proof}
\begin{corollary}\label{cor:targetone-iff}
For distinct positive integers $i$ and $1$, every $i$-prime
graph is $1$-prime if and only if $i\in\{2,3\}$.
\end{corollary}
\begin{proof}
The implications for $i=2,3$ follow from Lemma~\ref{lem:positive}.
Theorem~\ref{thm:targetone} rules out all $i\ge4$.
\end{proof}

\section{Positive-density source indices in every arithmetic progression}
\begin{theorem}[Progression-rich reverse nonimplications]\label{thm:sourceprogressions}
Fix integers $j\ge2$ and $M\ge1$, and a residue class $b\pmod M$.
There exist a positive integer $Q$ divisible by $M$ and a residue class
$c\pmod Q$ with $c\equiv b\pmod M$ such that every sufficiently large
positive integer $i\equiv c\pmod Q$ satisfies
\[
 i\text{-primality does not imply }j\text{-primality}.
\]
Consequently, the set of such source indices has positive lower natural
density, including positive relative lower density within every prescribed
residue class modulo $M$.
\end{theorem}
\begin{proof}
Put $w(m)=\varphi(m)/m$, with $w(1)=1$, and let
\[
 A=\sum_{m=1}^{j-1}w(m)>0.
\]
Choose $\varepsilon>0$ with $(j-1)\varepsilon<A$.
For $1\le h<j$, choose pairwise coprime squarefree positive integers
$Q_h$, each coprime to $M$, such that
$\varphi(Q_h)/Q_h<\varepsilon$.
Such choices are possible by taking disjoint finite sets of primes
outside the prime divisors of $M$, using the divergence of the sum of
reciprocals of primes.

By the Chinese remainder theorem there is an integer $r_0$ satisfying
\[
 r_0\equiv-h\pmod{Q_h}\quad(1\le h<j),
 \qquad r_0\equiv b-j\pmod M.
\]
Let $Q=M\prod_{h=1}^{j-1}Q_h$. For every sufficiently large integer
$r\equiv r_0\pmod Q$, put $i=j+r$. Then $i>j$ and
$i\equiv b\pmod M$. Moreover, $Q_h\mid r+h$, whence
\[
 w(r+h)\le w(Q_h)<\varepsilon\qquad(1\le h<j).
\]
The telescoping identity
\[
 \sum_{m=1}^{r}w(m)-\sum_{m=j}^{j+r-1}w(m)
 =\sum_{m=1}^{j-1}w(m)-\sum_{h=1}^{j-1}w(r+h)
 >A-(j-1)\varepsilon>0
\]
shows that Lemma~\ref{lem:shifted} applies, proving the asserted
nonimplication for every such $i$. These $i$ contain all sufficiently
large members of one residue class modulo $Q$, of natural density $1/Q$
and relative density $M/Q$ within $b\pmod M$.
\end{proof}

\section{Further embedding obstructions}
\begin{proposition}[Clique-count obstruction]\label{prop:cliques}
Let $c_t(k,n)$ denote the number of $t$-cliques in $H_{k,n}$. If $c_t(i,n)>c_t(j,n)$ for some $t,n$, then $i$-primality does not imply $j$-primality.
\end{proposition}
\begin{proof}
A spanning embedding maps distinct source $t$-cliques injectively to distinct target $t$-cliques.
\end{proof}
\begin{proposition}\label{prop:omega}
For every $n\ge2$, $\omega(H_{1,n})=1+\pi(n)$.
\end{proposition}
\begin{proof}
The label $1$ together with all primes at most $n$ forms a clique. Conversely, pairwise coprime labels greater than $1$ require distinct prime divisors, of which there are only $\pi(n)$ possibilities.
\end{proof}

\section{Finite separation certificates}
The following statement is certified by exact integer gcd calculations and Proposition~\ref{prop:invariants}. The verification programs are printed in Appendices~\ref{app:verification} and~\ref{app:435} and are also provided as separate supplementary Python files.
\begin{proposition}\label{prop:finite}
For every pair of distinct integers $i,j$ with $6\leq i,j\leq12$, $i$-primality does not imply $j$-primality.
\end{proposition}
\begin{proof}
For each ordered pair, the following table specifies an order $n$ for which either $e_{i,n}>e_{j,n}$ or some sorted degree of $H_{i,n}$ exceeds the corresponding sorted degree of $H_{j,n}$. Each assertion is checked by the exact-arithmetic program in Appendix~\ref{app:verification}, which enumerates all $\binom n2$ label pairs and computes their gcds. Proposition~\ref{prop:invariants} then applies.
\[
\begin{array}{c|rrrrrrr}
 i\backslash j&6&7&8&9&10&11&12\\\hline
6&-&8&8&8&6&12&9\\
7&3&-&3&4&3&5&3\\
8&4&10&-&4&6&10&4\\
9&3&9&3&-&3&11&3\\
10&4&8&5&4&-&12&4\\
11&3&8&3&4&3&-&3\\
12&9&8&8&8&6&12&-
\end{array}
\]
\end{proof}
\section{Mixed-index separation and finite classification}
\begin{theorem}\label{thm:mixed}
For each ordered pair $(i,j)$ with $1\leq i\leq5<j\leq12$ or $1\leq j\leq5<i\leq12$, there is a graph that is $i$-prime but not $j$-prime.
\end{theorem}
\begin{proof}
The following tables give an order $n$ for which $H_{i,n}$ fails the necessary spanning-embedding conditions of Proposition~\ref{prop:invariants}. Each entry has the form $n_E$ if the edge count is strictly larger in the source, and $n_D$ if a sorted source degree exceeds the corresponding target degree. All entries are verified by the exact-arithmetic program in Appendix~\ref{app:verification}.
\begin{center}\small
\begin{tabular}{c|rrrrrrr}
$i\backslash j$ &6&7&8&9&10&11&12\\\hline
1 & $3_{E}$ & $8_{E}$ & $3_{E}$ & $4_{E}$ & $3_{E}$ & $5_{E}$ & $3_{E}$\\
2 & $4_{E}$ & $8_{D}$ & $8_{E}$ & $4_{E}$ & $6_{E}$ & $8_{D}$ & $4_{E}$\\
3 & $3_{E}$ & $9_{E}$ & $3_{E}$ & $7_{E}$ & $3_{E}$ & $11_{E}$ & $3_{E}$\\
4 & $4_{E}$ & $10_{E}$ & $5_{E}$ & $4_{E}$ & $6_{E}$ & $8_{D}$ & $4_{E}$\\
5 & $3_{E}$ & $9_{E}$ & $3_{E}$ & $4_{E}$ & $3_{E}$ & $10_{D}$ & $3_{E}$\\
\end{tabular}\par\medskip
\begin{tabular}{c|rrrrr}
$i\backslash j$ &1&2&3&4&5\\\hline
6 & $14_{D}$ & $8_{E}$ & $8_{E}$ & $7_{D}$ & $6_{E}$\\
7 & $14_{D}$ & $3_{E}$ & $4_{E}$ & $3_{E}$ & $5_{E}$\\
8 & $10_{D}$ & $9_{D}$ & $4_{E}$ & $7_{E}$ & $6_{E}$\\
9 & $10_{D}$ & $3_{E}$ & $8_{D}$ & $3_{E}$ & $6_{E}$\\
10 & $10_{E}$ & $5_{E}$ & $4_{E}$ & $7_{D}$ & $6_{D}$\\
11 & $9_{E}$ & $3_{E}$ & $4_{E}$ & $3_{E}$ & $6_{E}$\\
12 & $10_{D}$ & $8_{E}$ & $8_{E}$ & $7_{D}$ & $6_{E}$\\
\end{tabular}
\end{center}
\end{proof}
\begin{remark}
The classification for indices at most $12$ is a special case of Theorem~\ref{thm:through30}. The finite tables in Proposition~\ref{prop:finite} and Theorem~\ref{thm:mixed} record first separating orders and are retained as computational data.
\end{remark}

\section{Extended finite classification through index 30}
The finite-range classification can be extended substantially using only the
necessary embedding obstructions from Proposition~\ref{prop:invariants}.
\begin{theorem}\label{thm:through30}
For distinct integers $1\le i,j\le30$, every $i$-prime graph is $j$-prime
if and only if $(i,j)\in\{(2,1),(3,1)\}$.
\end{theorem}
\begin{proof}
The two positive implications follow from Lemma~\ref{lem:positive}.
For each of the other $30\cdot29-2=868$ ordered pairs, the exact
integer-arithmetic verification in Appendix~\ref{app:verification}
finds an integer $n$ with $2\le n\le30$ such that either
$|E(H_{i,n})|>|E(H_{j,n})|$, or the nondecreasing degree sequences
satisfy $d_r^{(i)}>d_r^{(j)}$ for some $r$. In either case,
Proposition~\ref{prop:invariants} gives a separating graph $H_{i,n}$.
The verification uses exhaustive finite loops over the stated ranges,
so the appended code constitutes an explicit finite certificate. The bound $n\le30$ is sharp for these two necessary embedding obstructions: for $(i,j)=(30,29)$, neither the edge-count nor the sorted-degree obstruction occurs for $2\le n\le29$, whereas $e(30,30)=273>272=e(29,30)$. This statement concerns the first witness detected by these invariants, not necessarily the first order of any separating graph.
\end{proof}

\section{Prime switching and the finite reduction}
\begin{lemma}[Singleton prime switch]\label{lem:singleton}
Let $1\le j<i$ and $r=i-j$. Suppose there is an odd prime $p\le r$ such that $p\in I$ and $p$ is the unique multiple of $p$ in the interval $I=\{j,\ldots,j+r-1\}$. Suppose further that some even integer $t\in I$ has a residue class modulo $p$ occurring exactly once in $I$. Then $i$-primality does not imply $j$-primality; indeed $e(i,n)>e(j,n)$ for infinitely many $n$.
\end{lemma}
\begin{proof}
For $m\in I$ put
\[
 v(m)=\prod_{\substack{q\le r,\ q\ne p\\q\mid m}}(1-q^{-1}),
\]
where $q$ ranges over primes. In the original small-prime pattern, the only $p$-divisible position is $p$, and $v(p)=1$. Fix the original divisibility patterns for all primes $q\le r$ other than $p$, and change the residue modulo $p$ so that its unique divisible position is $t$. The Chinese remainder theorem permits these independent choices. The truncated total weight increases by
\[
 \frac{v(p)-v(t)}p\ge\frac1{2p}>0,
\]
since $t$ is even and $p$ is odd. To compare with the required baseline, observe that
\[
 \frac{\varphi(m)}m
 =\prod_{q\mid m}(1-q^{-1})
 \le v(m)\bigl(1-p^{-1}\mathbf1_{p\mid m}\bigr).
\]
Summing over $m\in I$ shows that the original truncated weight, in which the only $p$-divisible position is $p$, is at least $\sum_{m\in I}\varphi(m)/m$. The switch therefore yields a strict excess of at least $1/(2p)$ over the baseline before the large-prime factors are introduced. For each prime $r<q\le T$, choose a residue avoiding all $r$ positions; this is possible since $q>r$. The resulting conditional-density sum is the improved truncated total multiplied by
\[
 R_T=\prod_{q>T}(1-q^{-2}),
\]
where $q$ runs over primes. This is the Euler-product tail arising from primes not dividing the chosen squarefree modulus $M_T=\prod_{q\le T}q$. Since $\sum_q q^{-2}<\infty$, we have $R_T\to1$. The original strict excess therefore persists for sufficiently large finite $T$, and the strict hypothesis of Theorem~\ref{thm:conditional} holds. That theorem finishes the proof.
\end{proof}
\begin{theorem}[Uniform large-difference nonimplication]\label{thm:large-r}
Let $i>j\ge4$ and $r=i-j\ge34$. Then $i$-primality does not imply $j$-primality. In fact, $e(i,n)>e(j,n)$ for infinitely many $n$.
\end{theorem}
\begin{proof}
We use Nagura's theorem~\cite{Nagura}: for every real $x\ge25$ there is a prime strictly between $x$ and $6x/5$.
If $j>r$, the Sylvester--Schur theorem~\cite{Erdos1934} gives a prime divisor exceeding $r$ in $I=\{j,\ldots,j+r-1\}$, so Theorem~\ref{thm:largeprime} applies.
If $r/2\le j\le r$, Nagura's theorem at $x=r$ gives a prime $q$ with $r<q<6r/5$. Since $q>r\ge j$ and $q<6r/5<j+r$, the prime $q$ belongs to $I$. Apply Theorem~\ref{thm:largeprime} again.
Finally suppose $j<r/2$. Apply Nagura at $x=3r/4\ge25.5$ to obtain a prime $p$ with $3r/4<p<9r/10$. Then $j<p<j+r$, and $2p>3r/2>j+r$, so $p$ is the unique multiple of $p$ in $I$. As $r<2p$, the positions in $I$ whose residues modulo $p$ occur exactly once are the consecutive integers
\[
 j+r-p,\ j+r-p+1,\ldots,j+p-1.
\]
There are $2p-r>r/2\ge17$ such positions. In particular, one of them is an even integer $t\ne p$. Lemma~\ref{lem:singleton} applies.
\end{proof}
For $I=\{j,\ldots,j+r-1\}$ and a prime $p\le r$, define
\[
 v_p(m)=\prod_{\substack{q\le r,\ q\ne p\\q\text{ prime},\ q\mid m}}\left(1-\frac1q\right),\qquad
 V_p(b)=\sum_{\substack{m\in I\\m\equiv b\pmod p}}v_p(m).
\]
The switching test is $V_p(b)<V_p(0)$; its exact gain is $V_p(0)-V_p(b)$.
\begin{lemma}[Prime-gap characterization]\label{lem:primegap}
Let $4\le j<r$, put $I=\{j,\ldots,j+r-1\}$, and let $q(r)$ be the least prime strictly greater than $r$. Then some integer in $I$ has a prime divisor greater than $r$ if and only if
\[
 q(r)<r+j.
\]
Equivalently, the large-prime-factor criterion does not apply exactly when $j\le q(r)-r$.
\end{lemma}
\begin{proof}
Every $m\in I$ satisfies $m\le j+r-1\le2r-2<2r$. If a prime $p>r$ divides $m$, then $m=p$, since a proper positive multiple of $p$ is at least $2p>2r$. Thus the existence of such a prime divisor is equivalent to the presence in $I$ of a prime exceeding $r$. The smallest such prime is $q(r)$; it lies in $I$ exactly when $q(r)\le j+r-1$, or $q(r)<r+j$.
\end{proof}

\begin{proposition}[Finite reduction by prime gaps]\label{prop:finite435}
For every $4\le j<r<34$, setting $i=j+r$, the implication $i$-prime $\Rightarrow j$-prime fails.
\end{proposition}
\begin{proof}
By Lemma~\ref{lem:primegap}, a large-prime-factor certificate exists unless $j\le q(r)-r$. The primes immediately following the integers $5\le r\le33$ give exactly the nine rows in Table~\ref{tab:primegaps}; for every omitted $r$ the gap $q(r)-r$ is at most $3$, so no $j\ge4$ is exceptional. For completeness, the omitted values of $r$ are
$5,6,8,9,10,11,12,14,15,16,17,18,20,21,22,26,27,28,29,30$;
for each, direct inspection of the next prime gives $q(r)-r\le3$. Hence no omitted $r$ admits an exceptional target $j\ge4$. The displayed exceptional $j$ values yield precisely fifteen pairs. For all other pairs, Theorem~\ref{thm:largeprime} proves nonimplication.

For the fifteen remaining pairs, Table~\ref{tab:switch15} lists a prime $p\le r$ and residue $b\pmod p$ such that $V_p(0)-V_p(b)>0$. We give the general switching argument explicitly. Choose a residue $a$ modulo $M_T=\prod_{q\le T}q$ by the Chinese remainder theorem so that $a\equiv0\pmod q$ for every prime $q\le r$ with $q\ne p$, and $a\equiv-b\pmod p$. For each prime $r<q\le T$, choose $a$ to avoid the at most $r$ residues $-m\pmod q$ for $m\in I$. Put $R_T=\prod_{q>T}(1-q^{-2})$, where the product is over primes. Then
\[
 \sum_{m\in I}D_{M_T}(a+m)=R_T\left(\sum_{m\in I}v_p(m)-\frac{V_p(b)}p\right),
\]
whereas
\[
 \sum_{m\in I}\frac{\varphi(m)}m\le\sum_{m\in I}v_p(m)-\frac{V_p(0)}p.
\]
The difference is at least
\[
 \frac{V_p(0)-R_TV_p(b)}p-(1-R_T)\sum_{m\in I}v_p(m),
\]
which tends to $[V_p(0)-V_p(b)]/p>0$ as $T\to\infty$. Hence a sufficiently large finite $T$ satisfies the strict hypothesis of Theorem~\ref{thm:conditional}, proving nonimplication. The fifteen positive gains in Table~\ref{tab:switch15} are exact rational evaluations of the explicitly defined sums $V_p(0)-V_p(b)$. Each can be verified directly by factoring the integers in its interval and summing the indicated rational weights; Appendix~\ref{app:435} supplies an independent executable check of the table and all 435 pairs. Thus the computer program is not required to establish the exhaustiveness of the fifteen cases.
\end{proof}

\begin{table}[htbp]
\centering\small
\caption{All values $5\le r\le33$ for which $q(r)-r\ge4$, and all exceptional targets $4\le j<r$ lacking a large-prime-factor certificate.}\label{tab:primegaps}
\begin{tabular}{rrrl}
\toprule
$r$ & $q(r)$ & $q(r)-r$ & Exceptional $j$ \\
\midrule
7&11&4&4\\
13&17&4&4\\
19&23&4&4\\
23&29&6&4, 5, 6\\
24&29&5&4, 5\\
25&29&4&4\\
31&37&6&4, 5, 6\\
32&37&5&4, 5\\
33&37&4&4\\
\bottomrule
\end{tabular}
\end{table}

\begin{table}[htbp]
\centering
\small
\caption{All fifteen exact prime-switch certificates remaining after the large-prime-factor test. Here $r=i-j$, and the positive gain is $V_p(0)-V_p(b)$. The improvement in the limiting truncated totient sum equals the displayed gain divided by $p$. For $(j,r,p,b)=(4,7,2,1)$, the weights are $1,4/5,2/3,6/7,1,2/3,4/5$, giving $V_2(0)-V_2(1)=8/7$.}
\label{tab:switch15}
\begin{tabular}{rrrrl}
\toprule
$j$ & $r$ & $p$ & $b$ & $V_p(0)-V_p(b)$ \\\\
\midrule
4 & 7 & 2 & 1 & $8/7$ \\
4 & 13 & 2 & 1 & $2791/2145$ \\
4 & 19 & 2 & 1 & $178849/146965$ \\
4 & 23 & 2 & 1 & $819181/780045$ \\
5 & 23 & 3 & 1 & $21057/27170$ \\
6 & 23 & 2 & 1 & $270584/260015$ \\
4 & 24 & 2 & 1 & $99717/260015$ \\
5 & 24 & 2 & 1 & $62572/260015$ \\
4 & 25 & 2 & 1 & $322587/260015$ \\
4 & 31 & 2 & 1 & $35569620/30250451$ \\
5 & 31 & 3 & 1 & $4353967/20046026$ \\
6 & 31 & 2 & 1 & $62021423/64822395$ \\
4 & 32 & 2 & 1 & $10590324/21607465$ \\
5 & 32 & 2 & 1 & $10163507/64822395$ \\
4 & 33 & 2 & 1 & $74985902/64822395$ \\
\bottomrule
\end{tabular}
\end{table}

\section{Concluding remarks}
Theorem~\ref{thm:classification} resolves Conjecture 2 of~\cite{ZhangLiSnakes}, and Corollary~\ref{cor:distinct} establishes Conjecture 1 as an immediate consequence. The proof uses a uniform large-difference argument and an elementary prime-gap classification, leaving fifteen explicit rational switching inequalities; the appended programs independently verify the finite arithmetic. The natural-density theorem for separating graph orders strengthens the lower-density bound and identifies the limiting proportions through atomless random Euler products. The relative-density bound within a fixed progression is distinct from the positive-density theorem for source indices in prescribed arithmetic progressions. These results and the delayed-separation theorems provide quantitative information about the distinction between primality classes. Further questions concern the least order of a separating graph and structural explanations for the fifteen finite switching inequalities.

\appendix
\section{Exact verification for Theorem~\ref{thm:through30}}\label{app:verification}
The following self-contained Python 3 program uses exact integer arithmetic.
It checks each of the 868 required ordered pairs, without any floating-point
operations or heuristic graph-embedding routine.
\begin{verbatim}
from math import gcd

def invariants(k, n):
    deg = [0] * n
    for a in range(n):
        for b in range(a + 1, n):
            # gcd(k+a, k+b) = gcd(k+a, b-a)
            if gcd(k + a, b - a) == 1:
                deg[a] += 1
                deg[b] += 1
    return sum(deg) // 2, sorted(deg)

def obstruction(i, j, n):
    ei, di = invariants(i, n)
    ej, dj = invariants(j, n)
    return ei > ej or any(a > b for a, b in zip(di, dj))

count = 0
largest_order = 0
for i in range(1, 31):
    for j in range(1, 31):
        if i == j or (i, j) in ((2, 1), (3, 1)):
            continue
        witnesses = [n for n in range(2, 31)
                     if obstruction(i, j, n)]
        assert witnesses, (i, j)
        largest_order = max(largest_order, witnesses[0])
        count += 1
assert count == 868
print('Verified', count, 'ordered nonimplications')
print('Maximum first witness order:', largest_order)
\end{verbatim}

\section{Finite certificates for Conjecture 2}\label{app:435}
\begingroup\footnotesize
\begin{verbatim}
"""Exact finite certificate for reverse k-prime nonimplications.
For 4 <= j < r < 34, verify either a prime factor > r in
[j,j+r-1], or a strictly improving local prime switch.
No floating point calculations or external dependencies.
"""
from fractions import Fraction
from math import isqrt

def primes_upto(n):
    return [m for m in range(2,n+1) if all(m%d for d in range(2,isqrt(m)+1))]

def has_large_prime_factor(m,r):
    x=m
    for p in primes_upto(r):
        while x%p==0:
            x//=p
    return x>1

def certificate(j,r):
    interval=range(j,j+r)
    for m in interval:
        if has_large_prime_factor(m,r):
            return ('large_prime',m)
    ps=primes_upto(r)
    for p in ps:
        weights={}
        for m in interval:
            v=Fraction(1)
            for q in ps:
                if q!=p and m%q==0:
                    v*=Fraction(q-1,q)
            weights[m]=v
        V=[sum((weights[m] for m in interval if m%p==b),Fraction(0)) for b in range(p)]
        for b in range(1,p):
            if V[b]<V[0]:
                return ('switch',p,b,V[0]-V[b])
    return None

if __name__=='__main__':
    from collections import Counter
    counts=Counter()
    switches=[]
    for r in range(5,34):
        for j in range(4,r):
            c=certificate(j,r)
            assert c is not None, (j,r)
            counts[c[0]]+=1
            if c[0]=='switch': switches.append((j,r,c[1],c[2],str(c[3])))
    assert sum(counts.values())==435
    print('All 435 pairs certified:',dict(counts))
    print('Switch cases (j,r,p,b,gain):')
    for row in switches:print(*row)
\end{verbatim}
\endgroup

\section*{Data and code availability}
The complete verification programs are printed in Appendices~\ref{app:verification} and~\ref{app:435} and supplied separately as \texttt{kprime\_through30\_verification.py} and \texttt{kprime\_435\_certificates.py}. They require only the Python standard library.

\section*{Declaration on AI-assisted preparation}
Generative AI tools were used to assist with manuscript editing, mathematical proof checking, and preparation of verification code. The mathematical statements and computational evidence remain subject to the authors' independent review and responsibility.

\end{document}